\documentclass [10pt,twoside,b5paper]{article}
\usepackage{amsfonts}
\usepackage{amsthm}
\usepackage{amsmath}
\usepackage{amstext}
\usepackage{amssymb}
\usepackage{mathrsfs}
\usepackage{epsf}
\usepackage{graphicx}
\usepackage{fancybox}
\usepackage{color}
\usepackage{fancyhdr}
\usepackage[hang,footnotesize]{caption2}

\def\mathcal{\mathscr}
\newfont{\aaa}{cmb10 at 18pt}
\newfont{\bbb}{cmb10 at 10pt}

\newcommand{\beq}{\begin{equation}}
\newcommand{\eeq}{\end{equation}}
\newcommand{\bey}{\begin{eqnarray}}
\newcommand{\eey}{\end{eqnarray}}
\newcommand{\beyy}{\begin{eqnarray*}}
\newcommand{\eeyy}{\end{eqnarray*}}

\begin{document}

%%%以下正文开始
\setcounter{page}{1}
\qquad\\[8mm]

\begin{center}
{\bf\Large Stochastic Partial Differential {\textcolor{red}{Equations}}\\
[2mm] Driven by
 Fractional L\'{e}vy Noises
}
\footnote{Project Supported by
 National  Natural Science Foundation of China with
granted No. 11001051, 11371010, {\textcolor{red}{10971249}},
11301263. Email: lvxuebin2008@163.com(X. L\"{u}),
nan5lu8@netra.nju.edu.cn(W. Dai).}
\end{center}
%\\[2mm]
\begin{center}{\bbb Xuebin L\"{u}$^{1,2}$, Wanyang Dai $^{1}$}
\end{center}
%\\[-1mm]
\noindent\footnotesize{1 Department of Mathematics, Nanjing
University, Nanjing, P. R. China 210093
\\
2
 Department of Applied Mathematics, College of Science, Nanjing
University of Technology, Nanjing, P. R. China 210009
\\[6mm]
\normalsize\noindent{\bbb Abstract}\quad In this paper,
{\textcolor{red}{we investigate stochastic partial differential
equations driven by multi-parameter anisotropic fractional L\'{e}vy
noises, including the stochastic Poisson equation, the linear heat
equation, and the quasi-linear heat equation. Well-posedness of
these equations under the fractional noises will be addressed. The}
multi-parameter anisotropic fractional L\'{e}vy noise
{\textcolor{red}{is defined}} as the formal derivative of the
anisotropic fractional L\'{e}vy random {\textcolor{red}{field. In
doing so, there are two folds involved. First, we consider}} the
anisotropic fractional L\'{e}vy random field as the generalized
functional of the path of the pure jump L\'{e}vy
{\textcolor{red}{process. Second, we build} the Skorohod integration
with respect to the multi-parameter anisotropic fractional L\'{e}vy
noise by white noise approach. \vspace{0.3cm}

\noindent{\bbb Keywords}\quad  {\textcolor{red}{Stochastic partial
differential equation;}} White noise analysis;  Pure-jump L\'{e}vy
process; Generalized L\'{e}vy random field;  Anisotropic fractional
L\'{e}vy random field; Anisotropic fractional L\'{e}vy
{\textcolor{red}{noise}}\\
{\bbb MSC}(2000)\quad 60E07, 60G20, 60G51, 60G52, 60H40 \\[0.4cm]
\setcounter{equation}{0} \noindent

\noindent{\bbb{1\quad Introduction}}\\[1mm]
The study on fractional processes started from the fractional
 Brownian motion     introduced by Kolmogrov
  \cite{Kol1} and popularized by Mandelbrot and Van Ness \cite{Mandelbrot}.  The
  self-similarity and long-range dependence properties make the fractional Brownian
  motion
 suitable to model driving noises in different applications such as
hydrology. However, {\textcolor{red}{to capture the large jumps and
to}} model the higher variability phenomena, it is  natural to
consider {\textcolor{red}{more}} general long-range dependent
processes. {\textcolor{red}{Replacing a Brownian motion with a
L\'{e}vy process to define a} fractional process
{\textcolor{red}{becomes more and more popular (see, e.g,}}
\cite{HuangLC} \cite{HuangLP}
\cite{Hli}\cite{HuangLU}\cite{Lb}\cite{Marquardt}). Particularly,
the fractional L\'{e}vy
 process is long-range dependent and its one-dimensional
 distribution is infinitely divisible.
{\textcolor{red}{Furthermore, in} order to use the fractional
L\'{e}vy processes to model the higher variability phenomena, it is
{\textcolor{red}{imperative}} to investigate the stochastic calculus
{\textcolor{red}{based on}} fractional L\'{e}vy processes and
{\textcolor{red}{study the related}} stochastic differential
equation. {\textcolor{red}{Along the line, the works presented in}}
\cite{LUHuang} and \cite{Marquardt} {\textcolor{red}{are concerned
with}} the stochastic integral for deterministic integrands with
respect to fractional L\'{e}vy processes; In
 \cite{Bender},  the authors investigate the Skorohod integral
for fractional L\'{e}vy process whose underlying L\'{e}vy process
has finite moment of any {\textcolor{red}{order; In~\cite{LuDai},
the authors study stochastic (ordinary) differential equations
driven by fractional L\'{e}vy noises.}}
\par
  Lokka and Proske
\cite{Lokka} developed {\textcolor{red}{the}} white noise calculus
for pure jump L\'{e}vy process by viewing it as an element in the
Poisson space. According to \cite{Lokka}, on the Poisson space every
square integrable functional of the path of a pure jump L\'{e}vy
process has a chaos expansion in terms of the Charlier polynomials
with respect to Poisson random {\textcolor{red}{measure. Thus, a
so-called}} S-transformation can {\textcolor{red}{be used to}}
characterize a generalized stochastic distribution.
{\textcolor{red}{Moreover,} by the definition of the multi-parameter
fractional L\'{e}vy random {\textcolor{red}{field,}} we can
{\textcolor{red}{consider}} it as a square integrable functional of
the path of a pure jump L\'{e}vy process.
{\textcolor{red}{Therefore,}} the multi-parameter fractional
L\'{e}vy random field has a chaos expansion {\textcolor{red}{, which
implies that we} can use the white noise calculus of the pure jump
L\'{e}vy process given by \cite{Lokka} to handle its stochastic
integration. {\textcolor{red}{Well-posedness of these equations will
be addressed.}}

\par
 Motivated by  the
white noise  analysis for pure jump L\'{e}vy process given by Lokka
and Proske \cite{Lokka}, we define the stochastic integration with
respect to the anisotropic fractional L\'{e}vy random
{\textcolor{red}{fields. Furthermore, based on the integration,}} we
investigate several kinds of stochastic partial differential
equations driven by anisotropic fractional L\'{e}vy noises including
Poisson equation, linear heat equation and quasi-linear heat
equation.
\par
 This paper is organized as follows: In Section 2, we
recall the basic results about the white noise analysis of the
square integrable pure jump L\'{e}vy process given by A. Lokka and
F. N. Proske \cite{Lokka}; Based on S-transformation,  the
multi-parameter fractional L\'{e}vy noises are  introduced in
Section 3 as  the formal derivative of anisotropic fractional
L\'{e}vy random fields, and the Skorohod integral with respect to
multi-parameter fractional L\'{e}vy noise is built.  After the
preparation,    we investigate stochastic Poisson
{\textcolor{red}{equation}} driven by
{\textcolor{red}{$d$}}-parameter fractional L\'{e}vy noises in
Section 4; In Section 5, we investigate stochastic linear heat
equation  driven by anisotropic fractional L\'{e}vy noises. In
Section 6, under Lipschtz and linear conditions we obtain a unique
solution for stochastic quasi-linear heat equation  driven by
anisotropic fractional L\'{e}vy {\textcolor{red}{noises.}}\\[4mm]
\noindent
\noindent{\bbb 2\quad White noise calculus for pure jump L\'{e}vy
process}\\[1mm]

\noindent {\textcolor{red}{For convenience to readers and
citation}}, in this section, we recall {\textcolor{red}{some}} basic
results of white noise analysis of the pure jump L\'{e}vy process on
the Poisson space given by Lokka and {\textcolor{red}{Proske
\cite{Lokka}.}}
 \par
{\textcolor{red}{First,}} we recall the construction of the Poisson
space which has the similar properties as the classical Schwartz
space. Let $\xi_{n}$ denote the $n'$th Hermite function, the set of
Hermite functions $\{\xi_{n}\}_{n\in\mathbb{N}}$ is an orthonormal
basis of $L^{2}(\mathbb{R})$. Denote by
$\mathcal{S}$$(\mathbb{R}^{d})$ the Schwartz space of rapidly
decreasing $C^{\infty}$-functions on $\mathbb{R}^{d}$ and by
$\mathcal{S}'$$(\mathbb{R}^{d})$ the space of tempered
distributions. The nuclear topology on
$\mathcal{S}$$(\mathbb{R}^{d})$ is induced by the pre-Hilbertian
norms
$$
\|\phi\|_{p}^{2}:=\sum_{\alpha=(\alpha_{1},
\ldots,\alpha_{d})\in\mathbb{N}^{d}
  }(1+\alpha)^{2p}(\phi,\xi_{\alpha})_{L^{2}(\mathbb{R}^{d})}^{2},p\in\mathbb{N}
  _{0},
$$
where $(1+\alpha)^{2p}=\prod_{i=1}^{d}(1+\alpha_{i})^{2p}$,
$\xi_{\alpha}(x_{1},
\ldots,x_{d})=\prod_{i=1}^{d}\xi_{\alpha_{i}}(x_{i}) $, $\mathbb{N}
  _{0}=\mathbb{N}\setminus
  \{0\}$.
 Let $\mathbb{U}=\mathbb{R}^{d}\times \mathbb{R}_{0}$,where $\mathbb{R}
  _{0}=\mathbb{R}\setminus
  \{0\}$,  define
$$ \mathcal{S}(\mathbb{U}):=\{\phi\in
\mathcal{S}(\mathbb{R}^{d+1}):\phi(x_{1},
\ldots,x_{d},0)=\frac{\partial\phi}{\partial x_{d+1}}(x_{1},
\ldots,x_{d},0)=0 \}.
$$
$\mathcal{S}(\mathbb{U})$ is a closed subspace of
$\mathcal{S}(\mathbb{R}^{ d+1})$, thus it is a countably Hilbertian
nuclear algebra endowed with the topology induced by the norms  $\|
\cdot\|_{p}$, and its dual $\mathcal{S}'(\mathbb{U})\supset
\mathcal{S}'(\mathbb{R}^{d+1})$. For
$\phi\in\mathcal{S}(\mathbb{U})$, $\Phi\in\mathcal{S}'(\mathbb{U})$,
the action of $\Phi$ on $\phi$ is given by
 $
 \langle\Phi,\phi\rangle=\int_{\mathbb{U}}\Phi(x)\phi(x)d\lambda^{\times
 (d+1)}(x),
 $
$\lambda^{\times d}$ is the Lebesgue measure on $\mathbb{R}^{d }$.
Assume that $\nu$ is the L\'{e}vy
 measure  on $\mathbb{R}_{0}$ satisfying
$$
 \int_{\mathbb{R}_{0}}|x|^{2}d\nu(x)<\infty.\eqno(2.1)
$$
Denote   $\pi$   the measure on $\mathbb{U}$ given by
$\pi=\lambda^{\times d}\times\nu$. By Lemma 2.1 of \cite{Lokka},
there exists an element denoted by $1\otimes\dot{\nu}$ in
$\mathcal{S}'(\mathbb{U})$ such that
$$
\langle 1\otimes\dot{\nu},
\phi\rangle=\int_{\mathbb{U}}\phi(x)\pi(dx),\phi\in\mathcal{S}(\mathbb{U}).\eqno(2.2)
$$
In a generalized sense,  $1\otimes\dot{\nu}$ is the Randon-Nikodym
derivative of $\pi$ with respect to the Lebesgue measure.
\par
Denote ${L}^{2}(\mathbb{U} ,\pi )$ by the space of all square
integrable functions on $\mathbb{U} $ with respect to $\pi $, let
$(\cdot,\cdot)_{\pi}$ be the inner product on
${L}^{2}(\mathbb{U},\pi)$ and $|\cdot|_{\pi}$ the corresponding
norms on this space.
\par
 Define  $
\mathcal{N}_{\pi}:=\{\phi\in\mathcal{S}(\mathbb{U}):|\phi|_{\pi}=0\}
 $,  then $\mathcal{N}_{\pi}$ is a closed ideal of
$\mathcal{S}(\mathbb{U})$.
 Let $\widetilde{\mathcal{S}}(\mathbb{U})$ be
the space $
\widetilde{\mathcal{S}}(\mathbb{U}):=\mathcal{S}(\mathbb{U})/\mathcal{N}_{\pi}
$ endowed with the topology induced by the system of norms $
\|\widehat{\phi}\|_{p,\pi}:=\inf_{\psi\in\mathcal{N}_{\pi}}\|\phi+\psi\|_{p},
$ then $\widetilde{\mathcal{S}}(\mathbb{U})$ is a nuclear algebra.
Let $\widetilde{\mathcal{S}}'(\mathbb{U})$ be the dual of
$\widetilde{\mathcal{S}}(\mathbb{U})$, and for $p\in\mathbb{N}$, let
$\widetilde{\mathcal{S}}_{p}(\mathbb{U})$ denote the completion of
$\widetilde{\mathcal{S}}(\mathbb{U})$ with respect to the norm
$\|\cdot\|_{p,\pi}$, $\widetilde{\mathcal{S}}'_{-p}(\mathbb{U})$
denote the dual of $\widetilde{\mathcal{S}}_{p}(\mathbb{U})$.
$\widetilde{\mathcal{S}}(\mathbb{U})$ is the projective limit of
$\{\widetilde{\mathcal{S}}_{p}(\mathbb{U}), p> 0\}$, and
$\widetilde{\mathcal{S}}'(\mathbb{U})$ is the inductive limit of
$\{\widetilde{\mathcal{S}}'_{-p}(\mathbb{U}), p>0\}$.
$\widetilde{\mathcal{S}}(\mathbb{U})$ has similar nice properties as
the classical Schwartz space, so Lokka and Proske \cite{Lokka}
introduced it as the probability space to   construct   the white
noise analysis for  pure jump L\'{e}vy
process. \\
  \textbf{Theorem 2.1} (Lokka and Proske \cite{Lokka})\quad(1)There exists a probability measure $\mu_{\pi}$ on
$\widetilde{\mathcal{S}}'(\mathbb{U})$ such that
$$
\int_{\widetilde{\mathcal{S}}'(\mathbb{U})}e^{i\langle \omega,
\phi\rangle}d\mu_{\pi}(\omega)=\exp\{\int_{\mathbb{U}}(e^{i\phi(x)}-1)d\pi(x)\},
\forall\phi\in\widetilde{\mathcal{S}}(\mathbb{U}).\eqno(2.3)
$$
(2) Moreover, there exists a $p_{0}\in\mathbb{N}$ such that
$1\otimes\dot{\nu}\in
\widetilde{\mathcal{S}}_{-p_{0}}'(\mathbb{U})$, and a natural number
$q_{0}>p_{0}$ such that the imbedding operator
$\widetilde{\mathcal{S}}_{q_{0}}'(\mathbb{U})\hookrightarrow\widetilde{\mathcal{S}}_{p_{0}}'(\mathbb{U})$
is Hilbert-Schimidt and
$\mu_{\pi}(\widetilde{\mathcal{S}}_{-q_{0}}'(\mathbb{U}))=1$.
\par
From now on, for all $q_{0}$, $p_{0}$ are described in the Theorem
2.1. Set $\Omega=\widetilde{\mathcal{S}}'(\mathbb{U})$ and
$P=\mu_{\pi}$ given by Theorem 2.1,   Lokka and Proske \cite{Lokka}
give the infinite dimensional calculus for pure jump measure on
$(\Omega, P)$, and all of our following discussion is based on this
probability space.
\par
 Let $C_{n}(\cdot)$ be the
Charlier polynomials given by \cite{Lokka}, especially for $n=1$,
$C_{1}(\omega)=\omega- 1\otimes\dot{\nu},$  $\omega\in\Omega$. \\
\textbf{Lemma 2.2}(Lokka and Proske \cite{Lokka})\quad  For all $m,
n\in \mathbb{N}$,
$\varphi^{(n)}\in\widetilde{\mathcal{S}}(\mathbb{U})^{\widehat{\otimes}n}$,
$\psi^{(m)}\in
\widetilde{\mathcal{S}}(\mathbb{U})^{\widehat{\otimes}m}$,
($\widehat{\otimes}$  denotes the symmetrized tensor product),  the
following orthogonality relation holds,
\begin{equation}
\int_{\widetilde{\mathcal{S}}'(\mathbb{U})}\langle C_{n}(\omega),
\varphi^{(n)}\rangle\langle
C_{m}(\omega), \psi^{(m)}\rangle d\mu_{\pi}(\omega)=\begin{cases}0,&\  n\neq m\\
n!(\varphi^{(n)},\psi^{(n)})_{\pi},&\ n=m \ .\nonumber
\end{cases}
\end{equation}
\quad  Since $\widetilde{\mathcal{S}}(\mathbb{U})$  is dense in $L
^{2}(\mathbb{U})$, for $f \in L ^{2}(\mathbb{U})$, there exists a
sequence of functions $f _{n}\in \widetilde{\mathcal{S}}(\mathbb{U})
$ such that $f_{n}\rightarrow f $ in $L^{2}(\mathbb{U} ,\pi) $ as
$n\rightarrow\infty$. Define $\langle C_{1}(\omega), f\rangle$ by $
\langle C_{1}(\omega), f \rangle=\lim_{n\rightarrow\infty}\langle
C_{1}(\omega), f_{n}\rangle (limit \ in \ L^{2}(\mu_{\pi})), $ the
definition is independent of the choice of approximating sequence,
and by Lemma 2.2  the following isometry holds
$$
\int_{\widetilde{\mathcal{S}}'(\mathbb{U})}\langle
C_{1}(\omega),f\rangle^{2}
d\mu_{\pi}(\omega)=\int_{\widetilde{\mathcal{S}}'(\mathbb{U})}\langle
\omega-1\otimes\dot{\nu},f\rangle^{2}
d\mu_{\pi}(\omega)=|f|_{\pi}^{2}.\eqno(2.4)
$$
For any Borel sets $\Lambda_{1}\subset \mathbb{R}^{d}$ and
$\Lambda_{2}\subset \mathbb{R}_{0}$ such that the 0 is not in the
closure of $\Lambda_{2}$, define the random measure
$$
N(\Lambda_{1},\Lambda_{2}):=\langle
\omega,1_{\Lambda_{1}\times\Lambda_{2}}\rangle,
\widetilde{N}(\Lambda_{1},\Lambda_{2}):=\langle
\omega-1\otimes\dot{\nu},1_{\Lambda_{1}\times\Lambda_{2}}\rangle.
$$
From the characterization function of $\mu_{\pi}$, it is easy to
deduce that  $N$ is a Poisson random measure, and $\widetilde{N}$ is
the corresponding compensated measure. The compensator of $
N(\Lambda_{1},\Lambda_{2})$ is given by $\langle
1\otimes\dot{\nu},1_{\Lambda_{1}\times\Lambda_{2}}\rangle $ which is
equal to $\pi(\Lambda_{1}\times\Lambda_{2})$. Moreover,
$$
\int_{\mathbb{U}}\phi(s,x)\widetilde{N}(ds,dx)=\langle
\omega-1\otimes\dot{\nu},\phi\rangle, \phi\in
L^{2}(\mathbb{U},\pi).\eqno(2.5)
$$
By (2.2) and (2.3), we have
$$
\int_{\widetilde{\mathcal{S}}'(\mathbb{U})}e^{i\langle
\omega-1\otimes\dot{\nu},
\phi\rangle}d\mu_{\pi}(\omega)=\exp\{\int_{\mathbb{U}}(e^{i\phi(x)}-1-i\phi(x))d\pi(x)\},
\forall\phi\in\widetilde{\mathcal{S}}(\mathbb{U}).\eqno(2.6)
$$
 Denote $\mathcal{B}(\mathbb{R}^{d})$ the Borel
$\sigma$-algebra on $\mathbb{R}^{d}$, for $S\in
\mathcal{B}(\mathbb{R}^{d})$, define $X(S)$ by
$$
X(S)(\omega)=\langle C_{1}(\omega),\phi_{S}\rangle= \langle
\omega-1\otimes\dot{\nu},\phi_{S}\rangle,
$$
where
$\phi_{S}(x_{1},\ldots,x_{d},x_{d+1})=1_{_{S}}(x_{1},\ldots,x_{d})\times
x_{d+1} $. By (2.6), we have
$$\aligned
&\int_{\widetilde{\mathcal{S}}'(\mathbb{U})}e^{i X(S)}
d\mu_{\pi}(\omega)\\&
=\exp\{\int_{\mathbb{U}}(e^{i1_{_{S}}(x_{1},\ldots,x_{d})\times
x_{d+1}}-1-i1_{_{S}}(x_{1},\ldots,x_{d})\times x_{d+1})d\pi(x)\}\\
&= \exp\{\int_{S\times\mathbb{R}_{0}}(e^{i x_{d+1}}-1-i
x_{d+1})d\nu(x_{d+1})d\lambda^{d}(x_{1},\ldots,x_{d})\}\\
 &=\exp\{Leb(S)\int_{\mathbb{R}_{0}}(e^{iy}-1-iy)d\nu(y)\},
\endaligned \eqno(2.7)
$$
 where $Leb(S)$ is the Lebesgue measure of $S$, then $ \{X(S), S\in \mathcal{B}(\mathbb{R}^{d})\}$ is
a real-valued   random measure on $\mathbb{R}^{d}$. For $f\in
L^{2}(\mathbb{R}^{d})$, define
$$
\dot{X}(f)(\omega)=\langle C_{1}(\omega),\hat{f}\rangle= \langle
\omega-1\otimes\dot{\nu},\hat{f}\rangle,    \eqno(2.8)
$$
where $\hat
{f}(x_{1},\ldots,x_{d},x_{d+1})=f(x_{1},\ldots,x_{d})\times x_{d+1}
$.
 By (2.6), we obtain
$$
\int_{\widetilde{\mathcal{S}}'(\mathbb{U})}e^{i
\dot{X}(f)(\omega)}d\mu_{\pi}(\omega)=\exp\{ \int_{\mathbb{U}}
 (e^{i\hat{f}(x)}-1-i\hat{f}(x))d \pi(x)\}.\eqno(2.9)
$$
Hence,  we  can write formally
$$
\dot{X}(f)= \int_{\mathbb{R}^{d}}f(x)dX(x), f\in
L^{2}(\mathbb{R}^{d}).\eqno(2.10)
$$
Moreover, by (2.4),
$$
E_{\mu_{\pi}}(\dot{X}(f))^{2}=\|f\|_{L^{2}}^{2}\int_{\mathbb{R}_{0}}|x|^{2}d\nu(x).
\eqno(2.11)
$$
\par
 Now we recall the space of the stochastic distribution functions defined by \cite{Lokka}. Define the space $
\mathcal{P}(\widetilde{\mathcal{S}}'(\mathbb{U})) =
\{f:\widetilde{\mathcal{S}}'(\mathbb{U})\rightarrow \mathbb{C},
f(\omega)=\sum_{n=0}^{N}\langle \omega^{\otimes n},
\phi^{(n)}\rangle,$     $\omega\in\widetilde{
\mathcal{S}}'(\mathbb{U}), \phi^{(n)}\in\widetilde{
\mathcal{S}}(\mathbb{U})^{\widehat{\otimes}n}, N\in\mathbb{N}\}
 $,
 $ f$ is called a continuous polynomial function if
$f\in\mathcal{P}(\widetilde{\mathcal{S}}'(\mathbb{U}))$  and it
admit a unique representation of the form
$$
f(\omega)=\sum_{n=0}^{\infty}\langle C_{n}(\omega),
f_{n}\rangle,f_{n}\in\widetilde{
\mathcal{S}}(\mathbb{U})^{\widehat{\otimes}n}.
$$
For any number  $p\geq q_{0}$, define the Hilbert space
$(\mathcal{S})_{p}^{1}$ as the completion of
$\mathcal{P}(\widetilde{\mathcal{S}}'(\mathbb{U}))$ with respect to
the norm
$$
\|f\|^{2}_{p,1}=\sum_{n=0}^{\infty}(n!)^{2}\|f_{n}\|_{p,\pi}^{2}.
$$
The corresponding inner product is
$$
((f,g))_{p,1}=\sum_{n=0}^{\infty}(n!)^{2}((f_{n},g_{n}))_{p,\pi}.
 $$
where $((\cdot,\cdot))_{p,\pi}$ denote the inner product on
$\widetilde{\mathcal{S}}_{p}(\mathbb{U})^{\widehat{\otimes}n}$.
Obviously, $(\mathcal{S})_{p+1}^{1}\subset(\mathcal{S})_{p}^{1}$,
\cite{Lokka} define $(\mathcal{S})^{1}$ as the projective limit of
$\{(\mathcal{S})_{p}^{1}, p\geq q_{0}\}$, and  $(\mathcal{S})^{1}$
is a nuclear Fr\'{e}chet space which can be densely imbedding in
$L^{2}(\mu_{\pi})$.  Denote $(\mathcal{S})_{-p}^{-1}$ as the dual of
$(\mathcal{S})_{p}^{1}$, $(\mathcal{S})^{-1}$ as the inductive limit
of $\{(\mathcal{S})_{-p}^{-1}, p\geq q_{0}\}$ which is equal to the
dual of  $(\mathcal{S})^{1}$. $F\in(\mathcal{S})^{-1}$ if and only
if $F$ admit an expansion
$$
F(\omega)=\sum_{n=0}^{\infty}\langle C_{n}(\omega),
F_{n}\rangle,F_{n}\in\widetilde{
\mathcal{S}}'(\mathbb{U})^{\widehat{\otimes}n},
$$
and there exists a $ p\geq q_{0}$ such that
$$
\|F\|^{2}_{-p,-1}=\sum_{n=0}^{\infty}\|F_{n}\|_{-p,\pi}^{2}<\infty .
$$
For  $F\in(\mathcal{S})^{-1} $,  $f\in(\mathcal{S})^{1}$,
$$
\langle\langle F,f\rangle\rangle=\sum_{n=0}^{\infty}n!\langle
F_{n},f_{n}\rangle_{\pi},
$$
 $ \langle\langle
\cdot,\cdot\rangle\rangle$ is an extension of the inner product on
$L^{2}(\mu_{\pi})$. $(\mathcal{S})^{1}$ is called space of
stochastic test functions, $(\mathcal{S})^{-1}$ is called space of
stochastic distribution functions, they are pairs of dual spaces and
$(\mathcal{S})^{1}\subset
L^{2}(\mu_{\pi})\subset(\mathcal{S})^{-1}$.
\par
Next we recall the S-transform given by \cite{Lokka} which can
transform stochastic distribution functions to deterministic
functionals. Let
$$
\widetilde{e}(\phi,\omega):=\exp(\langle\omega,\ln(1+\phi)\rangle-\langle1\otimes\dot{\nu},\phi\rangle),\eqno(2.12)
$$
it is analytic as a function of
$\phi\in\widetilde{\mathcal{S}}_{q_{0}}$ satisfying $\phi(x)>-1$ for
all $x\in \mathbb{U}$. Moreover, it has the following chaos
expansion,
$$
\widetilde{e}(\phi,\omega)=\sum_{n=0}^{\infty}\frac{1}{n!}\langle
C_{n}(\omega), \phi^{\otimes n}\rangle.\eqno(2.13)
$$
  Denote
 $ U_{p}:=\{\phi\in\widetilde{\mathcal{S}}(\mathbb{U}):
\|\phi\|_{p,\pi}<1\}  $, by (2.13), \cite{Lokka} proved that
$\widetilde{e}(\phi,\omega)\in(\mathcal{S})^{1}_{p}$ if and only if
$\phi\in U_{p}$.\\
 \textbf{Definition 2.3}(\cite{Lokka})\ Let $F\in(\mathcal{S})^{-1}_{-p} $, $\xi\in
 U_{p}$, the  S-transform of $F$ is defined by
 $$
 S(F)(\xi):=\langle\langle
 F,\widetilde{e}(\xi,\omega)\rangle\rangle.
$$
 For example, if $F=\sum_{n=0}^{\infty}\langle
C_{n}(\omega),F_{n}\rangle\in(\mathcal{S})^{-1}_{-p} $, $\xi\in
 U_{p},$ then  $
 S(F)(\xi)=\sum_{n=0}^{\infty}\langle F_{n},  \xi^{\otimes
 n}\rangle_{\pi}.
 $ \par
 Denote $\mathcal{U}=Hol(0)$  the algebra of germs of functions
that are holomorphic in a neighborhood of 0. The S-transform is
isomorphic between $(\mathcal{S})^{-1} $ and $\mathcal{U}$.\\
 \textbf{Theorem 2.4}(\cite{Lokka}) \ If $F \in (\mathcal{S})^{-1}_{-p}$,
 then $S(F)\in\mathcal{U}$. Conversely, if $G\in\mathcal{U}$, there
 is  a uniquely defined distribution  $F \in (\mathcal{S})^{-1}_{-p}$
such that $G=S(F)$ on some neighborhood of $0$ in
$(\mathcal{S})^{-1}$.
\par
Since  $f, g\in\mathcal{U}$, then $fg\in\mathcal{U}$, then  by
Theorem 2.4, the following definition of  Wick product is well-defined.\\
 \textbf{Definition 2.5}(\cite{Lokka})
\ Let $F,G\in (\mathcal{S})^{-1}$, define the Wick product
$F\diamond G$ of $F$ and $G$ by
$$
F\diamond G=S^{-1}( S(F)S(G)).
$$
The Wick exponential of $F\in (\mathcal{S})^{-1}$ denoted by
$\exp^{\diamond}(F)$ is defined by
$$
\exp^{\diamond}(F):=\sum_{n=0}^{\infty}\frac{1}{n!}F^{\diamond
n},\eqno(2.14)
 $$
 whenever$ \sum_{n=0}^{\infty}\frac{1}{n!}F^{\diamond n}\in
 (\mathcal{S})^{-1}$. In this case,
  $$S(\exp^{\diamond} X)(\eta)=\exp[(S X)(\eta)].\eqno(2.15)$$
Last recall the Skorohod integral of pure jump processes given by
\cite{Lokka}. Let $F:\mathbb{U}\rightarrow (\mathcal{S})^{-1}$ be
the random fields with chaos expansion
$$
F(x)=\sum_{n=0}^{\infty}\langle C_{n}(\omega),F_{n}(\cdot,x)\rangle,
$$
where $
F_{n}(\cdot,x)\in\mathcal{S}'(\mathbb{U})^{\widehat{\otimes}n}$ and
$\|F(x)\|_{-p,-1}<\infty, $ for some $p>0$. Let $\mathbb{L}$ denote
the set of all $F:\mathbb{U}\rightarrow (\mathcal{S})^{-1}$ such
that $
\widetilde{F_{n}}\in\mathcal{S}'(\mathbb{U})^{\widehat{\otimes}(n+1)}$($
\widetilde{F_{n}}$ is the symmetrization of $F_{n}$) and $
\sum_{n=0}^{\infty}|\widetilde{F_{n}}|_{-p,\pi}^{2}<\infty  $
for some $p>0$.\\
 \textbf{Definition 2.6}(\cite{Lokka})(Skorohod integral ) For
 $F\in\mathbb{L}$, define the Skorohod integral $\delta(F)$ by
 $$
 \delta(F):=\sum_{n=0}^{\infty}\langle
 C_{n+1}(\omega),\widetilde{F_{n}}\rangle.
$$
\ \ From the assumption on $\mathbb{L}$,  we see that
$\delta(F)\in(\mathcal{S})^{-1}$. For the predictable integrands,
the Skorohod integral coincides with the usual Ito-type integral
with respect to the compensated Poisson random measure.  \\
 \textbf{Proposition 2.7}(\cite{Lokka}) If $F\in\mathbb{L}$, then
$\delta(F)\in(\mathcal{S})^{-1}_{-p}$ for some $p>0$ and
$$
S\delta(F)(\xi)=\int_{\mathbb{U}} S F(x)(\xi)\xi(x)\pi(dx), \xi\in
U_{p}.\eqno(2.16)
$$
\\[4mm]
\noindent {\bbb 3\quad Anisotropic   fractional L\'{e}vy noises  }\\[1mm]

\noindent In this section, we define the  anisotropic  fractional
L\'{e}vy random field which can be considered as a generalized
functional of the path of the pure jump L\'{e}vy random field and
according to the result of section 2 we give its S-transformation.
Moreover, based on the S-transformation of  the  anisotropic
fractional L\'{e}vy random field, we define  its formal derivative
as {\textcolor{red}{$d$-parameter}} fractional L\'{e}vy noise.
\par Let $\overline{\beta}=(\beta_{1}, \ldots, \beta_{d}), 0<\beta_{k}<\frac{1}{2}, k=1,2,\ldots, d$, $f\in
\mathcal{S}(\mathbb{R}^{d})$, $\Gamma(\overline{\beta})=
\prod_{k=1}^{d} \Gamma(\beta_{k}) $, $ x
^{\overline{\beta}}=\prod_{k=1}^{d}x_{k}^{ \beta_{k}},x=(x_{1},
\ldots, x_{d}),$ the multi-variate fractional integral operator of
Liouville-type is defined by Samko et al.\cite{fraop}:
$$
I^{\overline{\beta} }_
{+\ldots+}f(x):=\frac{1}{\Gamma(\overline{\beta}) }
\int_{\mathbb{R}^{d}_{+}}\frac{f(x- y)dy}{ y ^{1-\overline{\beta}}},
\eqno(3.1)
$$
$$
I^{\overline{\beta} }_
{-\ldots-}f(x):=\frac{1}{\Gamma(\overline{\beta}) }
\int_{\mathbb{R}^{d}_{+}}\frac{f(x+y)dy}{y ^{1-\overline{\beta}}}.
\eqno(3.2)
$$
\textbf{Theorem 3.1}(Samko et al.\cite{fraop})\quad The operator
$I^{\overline{\beta} }_ {\pm\ldots\pm}$
  is
 bounded from $L^{\overline{p}}(\mathbb{R}^{d})$ to $L^{\overline{q}}(\mathbb{R}^{d})$
 with $\overline{p}=(p_{1},\ldots,p_{d})$,
 $\overline{q}=(q_{1},\ldots,q_{d})$ if and only if
 $$
1<p_{k}<\frac{1}{\beta_{k}}, q_{k}=\frac{p_{k}}{1-\beta_{k}p_{k}},
k=1,2, \ldots,d,
$$
where $L^{\overline{p}}$ is the  Banach space of functions with
mixed norm
$$
\|f\|_{\overline{p}}=\{\int_{\mathbb{R}}\{\ldots
\{\int_{\mathbb{R}}[\int_{\mathbb{R}}|f(s_{1},\ldots,s_{d})|
^{p_{1}}ds_{1}]^{\frac{p_{2}}{p_{1}}}ds_{2}\}^{\frac{p_{3}}{p_{2}}}\ldots\}
^{\frac{p_{d}}{p_{d-1}}}ds_{d}\}^{\frac{1}{p_{d}}}<\infty.\eqno(3.3)
$$
Especially, for $p_{1}=\ldots=p_{d}=p$,
$L^{\overline{p}}(\mathbb{R}^{d})$ is equal to
$L^{p}(\mathbb{R}^{d})$.
\par
Take $q_{1}=\ldots=q_{d}=2,
p_{k}=\frac{1}{\frac{1}{2}+\beta_{k}}, k=1,2, \ldots,d,$ we deduce
that
 for  $\overline{\beta}=(\beta_{1},\ldots, \beta_{d})$,
$0<\beta_{k}<\frac{1}{2}$, $k=1,\ldots,d$,  the operator
$I^{\overline{\beta}}: \mathcal{S}(\mathbb{R}^{d})\rightarrow
 L^{2}(\mathbb{R}^{d})$ is continuous. Hence, by (2.8), we can define the generalized  anisotropic fractional L\'{e}vy random
 field as follows:  \\
\textbf{Theorem 3.2}\quad  For $\overline{\beta}=(\beta_{1},\ldots,
\beta_{d})$, $0<\beta_{k}<\frac{1}{2}$, $k=1,\ldots,d$,
$$
\dot{X}^{\overline{\beta}}(f):=\dot{X}(I^{\overline{\beta}}_{-\ldots-}f),
\ f \in \mathcal{S}(\mathbb{R}^{d})\eqno(3.4)
$$
is  a tempered real-valued generalized random field. We refer it as
the generalized  anisotropic fractional L\'{e}vy random field.\\
\textbf{Proof:} The proof is the same  to that of Theorem 3.3 of
\cite{Lb}, we omit here.
\par
In fact, by (2.8), for $f \in \mathcal{S}(\mathbb{R}^{d})$, $
\dot{X}^{\overline{\beta}}(f)$ can be represented as
$$
\dot{X}^{\overline{\beta}}(f)=\langle C_{1},
K^{\bar{\beta}}f\rangle,
$$
where
$$
(K^{\bar{\beta}}f)(x_{1}, \ldots, x_{d},
x_{d+1})=I^{\overline{\beta}}_{-\ldots-}f(x_{1}, \ldots,
x_{d})\times x_{d+1}.
$$
Moreover, by (2.9), we get
$$
\int_{\widetilde{\mathcal{S}}'(\mathbb{U})}e^{i
\dot{X}^{\overline{\beta}}(f)(\omega)}d\mu_{\pi}(\omega)=\exp\{
\int_{\mathbb{U}}
 (e^{i(K^{\bar{\beta}}f)(x)}-1-i(K^{\bar{\beta}}f)(x))d \pi(x)\}
.\eqno(3.5)
$$
Since $I^{\overline{\beta} }_ {-\ldots-}1_{[0,\overline{t}]}\in
L^{2}(\mathbb{R}^{d}), \overline{t}\in\mathbb{R}_{+}^{d}$, by (2.8)
we can define the anisotropic fractional L\'{e}vy random field as
follows:\\
 \textbf{Definition 3.3} \ The anisotropic fractional
L\'{e}vy random field is defined by
$$
X^{\overline{\beta}}_{\overline{t}}:=\dot{X}(I^{\overline{\beta}}_{-\ldots-}1_{[0,\overline{t}]}),
\overline{t}=(t_{1}, \ldots,t_{d}), t_{i}\geq0, i=1, 2, \ldots,
d\eqno(3.6)
$$
(3.6) can be represented as
$$
X^{\overline{\beta}}_{\overline{t}}=\int_{\mathbb{R}^{d}}I^{\overline{\beta}
}_ {-\ldots-}1_{[0,\overline{t}]}(\overline{s})dX(\overline{s})
=\int_{-\infty}^{t_{1}}\ldots\int_{-\infty}^{t_{d}}\prod_{k=1}^{d}[(t_{k}-s_{k})^{\beta_{k}}-(-s_{k})_{+}^{\beta_{k}}]dX(\overline{s}).
\eqno(3.7)$$
 From  (3.7) we see that the fractional integral
parameters along different time axis are different, thus the
fractional L\'{e}vy random field  $\{
X^{\overline{\beta}}_{\overline{t}},\overline{t}\in\mathbb{R}_{+}^{d}\}$
is anisotropic.
\par Since $I^{\overline{\beta} }_ {-\ldots-}1_{[0,\overline{t}]}\in
L^{2}(\mathbb{R}^{d}), \overline{t}\in\mathbb{R}_{+}^{d}$,  $
X_{t}^{\overline{\beta}}$ has the following representation
$$
X_{\overline{t}}^{\overline{\beta}}=\langle C_{1},
K^{\bar{\beta}}1_{[0,t]}\rangle=\delta(K^{\bar{\beta}}1_{[0,t]}).\eqno(3.8)
$$
Thus, by (2.16), we  get the S-transform of anisotropic fractional
L\'{e}vy random field
$$
SX_{\overline{t}}^{\overline{\beta}}(\eta)=\int_{\mathbb{R}^{d}}\int_{\mathbb{R}_{0}}y\eta(\overline{s},y
)I^{\overline{\beta} }_
{-\ldots-}1_{[0,\overline{t}]}(\overline{s})\nu(dy)d\overline{s},
\eta\in U_{p}, p> p_{0}. \eqno(3.9)
$$
On the other hand, by   the following fractional integral by parts
formula of operator $I_{\pm\ldots\pm}^{\bar{\beta}}$:
$$
\int_{\mathbb{R}^{d}}
f(\overline{s})I_{+\ldots+}^{\overline{\beta}}g(\overline{s})d\overline{s}=
\int_{\mathbb{R}^{d}}g(\overline{s})I_{-\ldots-}^{\overline{\beta}}f(\overline{s})d\overline{s},
\ f,g\in \mathcal{S}(\mathbb{R})\eqno(3.10)
$$
which can be extended to $f\in L^{\bar{p}}(\mathbb{R})$, $g\in
L^{\bar{r}}(\mathbb{R})$ with $p_{i}>1$, $r_{i}>1$ and
$\frac{1}{p_{i}}+\frac{1}{r^{i}}=1+\beta_{i}$, $i=1,\ldots,d$, (3.9)
can be written as
$$
\aligned SX^{\overline{\beta}}_{\overline{t}}(\eta)
&=\int_{\mathbb{R}^{d}}\int_{\mathbb{R}_{0}}1_{[0,\overline{t}]}(\overline{s})yI_{+\ldots+}^{\overline{\beta}}\eta(\cdot,y
)(\overline{s})\nu(dy)d\overline{s}\\
&=\int_{0}^{t_{1}}\ldots\int_{0}^{t_{d}}[\int_{\mathbb{R}_{0}}yI_{+\ldots+}^{\overline{\beta}}\eta(\cdot,y
)(\overline{s})\nu(dy)]d\overline{s}.\endaligned \eqno(3.11)
$$
Hence,
$$
\frac{\partial^{d}}{\partial t_{1}\ldots\partial
t_{d}}SX_{\overline{t}}^{\overline{\beta}}(\eta)=\int_{\mathbb{R}_{0}}yI_{+\ldots+}^{\overline{\beta}}\eta(\cdot,y)(\overline{t})\nu(dy),
\eta\in U_{p}, p> p_{0}. \eqno(3.12)
$$
We denote $\dot{X}^{\bar{\beta}}_{\overline{t}}$  the fractional
L\'{e}vy noise in the following sense:
$$
S\dot{X}^{\overline{\beta}}_{\overline{t}}(\eta)=\frac{\partial^{d}}{\partial
t_{1}\ldots\partial t_{d}}
SX^{\beta}_{\overline{t}}(\eta)=\int_{\mathbb{R}_{0}}yI_{+\ldots+}^{\overline{\beta}}\eta(\cdot,y)(\overline{t})\nu(dy),
\eta\in U_{p},  p> p_{0}, \bar{t}\in\mathbb{R}^{d}_{+}.\eqno(3.13)
$$
Next we prove that $\dot{X}^{\overline{\beta}}_{\overline{t}}$ is a
generalized stochastic distribution function and it has a chaos
representation:\\
 \textbf{Theorem
3.4}\ $\dot{X}^{\overline{\beta}}_{\overline{t}}\in
(\mathcal{S})^{-1}_{-p}$ for all $p>\max\{1,p_{0}\}$  and
$$
\dot{X}^{\overline{\beta}}_{\overline{t}}=\langle C_{1},
\lambda_{\overline{t}}\rangle, \eqno(3.14)
$$
where
$$\lambda_{\overline{t}}(\overline{u},y)=\frac{y(\overline{t}-\overline{u})_{+}^{\overline{\beta}-1}}{\Gamma(\overline{\beta})}
=\frac{y\prod_{k=1}^{d}(t_{k}-u_{k})_{+}^{\beta_{k}-1}}{\Gamma(\overline{\beta})}
$$
 \textbf{Proof:}  We first show that  $
\langle C_{1}, \lambda_{\overline{t}}\rangle\in
(\mathcal{S})^{-1}_{-p}$ for all $p>\max\{1,p_{0}\}$. By the
estimate
$$
\int_{\mathbb{R}}(t-u)_{+}^{\beta-1} \xi_{n}(u)du   \leq C
n^{\frac{2}{ 3}-\frac{\beta}{2}}.\eqno(3.15)
$$
from section 4 of \cite{Elliott} , where $C$ is a certain constant
independent of $t$,
$$ \aligned
\|\langle C_{1}, \lambda_{\overline{t}}\rangle\|_{-1,-p}^{2}
&=\frac{\int_{\mathbb{R}}|y|^{2}d\nu(y)}{\Gamma(\overline{\beta})}
\sum_{ \alpha=(\alpha_{1}, \ldots,\alpha_{d})\in\mathbb{ N}_{0}^{d}}
(\alpha+1)^{-2p}\langle(\overline{t}-\cdot)_{+}^{\beta-1},
\xi_{\alpha}\rangle_{L^{2}(\mathbb{R}^{d})}^{2}\\
&=A\sum_{ \alpha=(\alpha_{1}, \ldots,\alpha_{d})\in\mathbb{
N}_{0}^{d}}
\prod_{k=1}^{d}(\alpha_{k}+1)^{-2p}(\int_{\mathbb{R}}(t_{k}-u_{k})_{+}^{\beta_{k}-1} \xi_{\alpha_{k}}(u)du)^{2}\\
 &\leq AC\sum_{ \alpha=(\alpha_{1}, \ldots,\alpha_{d})\in\mathbb{ N}_{0}^{d}} \prod_{k=1}^{d}(\alpha_{k}+1)^{-2p+\frac{4}{3}-\beta_{k}}\\
&=A C\prod_{k=1}^{d}\sum_{ \alpha_{k}=1}^{\infty} (\alpha_{k}+1)^{-2p+\frac{4}{3}-\beta_{k}}\\
&< +\infty,  for\  p>\max\{1,p_{0}\},
\endaligned  \eqno(3.16)$$
where
$A=\frac{\int_{\mathbb{R}}|y|^{2}d\nu(y)}{\Gamma(\overline{\beta})}$
is a positive constant. Thus $ \langle C_{1},
\lambda_{\overline{t}}\rangle\in (\mathcal{S})^{-1}_{-p}$ for all
$p>\max\{1,p_{0}\}$. Next, we prove (3.14) holds. In fact,
$$
(I_{-\ldots-}^{\overline{\beta}}\delta_{\bar{t}})(\bar{s})=\frac{1}{\Gamma(\overline{\beta})
}
\int_{\mathbb{R}^{d}_{+}}\frac{\delta_{\bar{t}}(\bar{s}+\bar{u})d\bar{u}}{\bar{u}
^{1-\overline{\beta}}}=\frac{(\bar{t}-\bar{s})^{\overline{\beta}-1}}{\Gamma(\overline{\beta})}.\eqno(3.17)
$$
Taking S-transform of $ \langle C_{1},
\lambda_{\overline{t}}\rangle$,
$$\aligned
 S\langle C_{1}, \lambda_{\overline{t}}\rangle(\eta)
&=\int_{\mathbb{R}^{d}}\int_{\mathbb{R}_{0}}\frac{y(\overline{t}-\overline{s})_{+}^{\overline{\beta}-1}}{\Gamma(\overline{\beta})}
\eta(\overline{s},y)\nu(dy)d\bar{s}\\
&=\int_{\mathbb{R}^{d}}\int_{\mathbb{R}_{0}}yI_{-\ldots-}^{\overline{\beta}}\delta_{\bar{t}}(\bar{s})\eta(\overline{s},y)\nu(dy)d\bar{s}\\
&=\int_{\mathbb{R}^{d}}\int_{\mathbb{R}_{0}}y\delta_{\bar{t}}(\bar{s})I_{+\ldots+}^{\overline{\beta}}\eta(\cdot,y)(\overline{s})\nu(dy)d\bar{s}\\
&=\int_{\mathbb{R}_{0}}yI_{+\ldots+}^{\overline{\beta}}\eta(\cdot,y)(\overline{t})\nu(dy),
\eta\in U_{p},  p> p_{0}.\endaligned
$$
Hence, by (3.13), (3.14) holds when take S-transfrom by two
sides.\hfill$\Box$
 \par
 Now
 we define the  Skorohod integral for
$(\mathcal{S})^{-1}$-valued processes with respect to   $
X^{\overline{\beta}}$. First, we define
$(\mathcal{S})^{-1}$-valued integrals as follows:\\
 \textbf{Definition 3.5} Suppose  $F:\mathbb{R}_{+}^{d}\longrightarrow
 (\mathcal{S})^{-1}$ is a given function such that
 $
 \langle\langle F(x), f\rangle\rangle\in L^{1}(\mathbb{R}_{+}^{d},dx)$ for
 all $f\in(\mathcal{S})$, then $\int_{\mathbb{R}_{+}^{d}}F(x)dx$ is defined to
 be the unique element of $(\mathcal{S})^{-1}$ such that
$$
\langle\langle \int_{\mathbb{R}_{+}^{d}}F(x)dx,
f\rangle\rangle=\int_{\mathbb{R}_{+}^{d}}\langle\langle F(x),
f\rangle\rangle dx.\eqno(3.18)
$$
\textbf{Definition 3.6} \ Suppose that
$F:\mathbb{R}_{+}^{d}\longrightarrow
 (\mathcal{S})^{-1}$
 such that $F(\bar{s})\diamond
\dot{X}^{\bar{\beta}}_{\bar{s}} $ is $d\bar{s}-$ integrable in $
(\mathcal{S})^{-1}$. Then we define the Skorohod integral of $F$
with respect to $ X^{\bar{\beta}}$
  by
$$
\delta^{\bar{\beta}}(F):=\int_{\mathbb{R}_{+}^{d}}F(\bar{s})\delta
X^{\bar{\beta}}_{\bar{s}}:=\int_{\mathbb{R}_{+}^{d}}F(\bar{s})\diamond
\dot{X}^{\bar{\beta}}_{\bar{s}}d\bar{s}.\eqno(3.19)
$$
In particular, if $A\subset\mathbb{R}_{+}^{d}$ is a Borel set, then
$$
\int_{A}F(\bar{s})\delta
X^{\bar{\beta}}_{\bar{s}}:=\int_{\mathbb{R}_{+}^{d}}1_{A}(\bar{s})F(\bar{s})\diamond
\dot{X}^{\bar{\beta}}_{\bar{s}}d\bar{s}.\eqno(3.20)
$$
By definition 3.6 , we get\\
 \textbf{Proposition  3.7} \ Let $F: \mathbb{R}_{+}^{d}\longrightarrow (\mathcal{S})^{-1}$ be Skorohod integrable
with respect to $ X^{\beta}$, $Y\in(\mathcal{S})^{-1}$, then
$$
Y\diamond \delta^{\bar{\beta}}(F)=\delta^{\bar{\beta}}(Y\diamond
F),\eqno(3.21)
$$
the equation holds whenever one side exists.
\noindent  \\[4mm]
\noindent {\bbb 4\quad The stochastic Poisson equation driven by d-parameter fractional L\'{e}vy noise}\\[1mm]
\noindent In this section, we investigate the stochastic Poisson
equation driven by d-parameter fractional L\'{e}vy noise:
 $$
\begin{cases}
  \triangle U (x)= -\dot{X}^{\overline{\beta}}_{x}, x\in
  D,
\\
 U (x)=0, x\in\partial D.
\end{cases}\eqno(4.1)
$$
where $\triangle=\sum_{k=1}^{d}\frac{\partial^{^{2}}}{\partial
x_{k}^{2^{}}}$ is the Laplace operator in $\mathbb{R}^{d}$,
$D\subset \mathbb{R}_{+}^{d}$ is a given domain with regular
boundary and $\dot{X}^{\overline{\beta}}_{x}$ is the  d-parameter
fractional L\'{e}vy noise.\\
 \textbf{Theorem 4.1}\  The stochastic Poisson equation (4.1) has a
 unique continuous solution in $(\mathcal{S})^{-1}$.\\
\textbf{Proof}: Based on the corresponding solution in the deterministic case ( with $\dot{X}^{\overline{\beta}}_{x}$
replaced by a bounded deterministic function ), the solution of (4.1) will be
 $$
U(x)=\int_{D}G(x,y)\dot{X}^{\overline{\beta}}_{y}dy,\eqno(4.2)
$$
where G is the Dirichlet    Laplacian. We first prove that
$U(x)\in(\mathcal{S})^{-1}_{-p}$ for all $p>\max\{1,p_{0}\}$. By
(3.14),  (4.2) can be written as
$$
U(x)=\int_{D}G(x,y)\dot{X}^{\overline{\beta}}_{y}dy=\int_{D}G(x,y)\langle
C_{1}, \lambda_{y}\rangle dy.
$$
Then by (3.16) and the fact that $ G(x,\cdot)\in
L^{1}(\mathbb{R}^{d})$,  for $p>\max\{1,p_{0}\}$,
 $$
\|U(x)\|_{-1,-p}
 \leq \int_{D}\|\langle C_{1},\lambda_{y}\rangle\|_{-1,-p} |G(x,y)| dy <+\infty,
  $$
   that is $U(x)=\langle C_{1},
\int_{D}G(x,y)\lambda_{y}dy\rangle\in(\mathcal{S})^{-1}_{-p}$ for
all $p>\max\{1,p_{0}\}$ and the same estimate gives that $
U(x):\overline{D}\rightarrow (S)^{-1}$ is continuous. It is easy to
show that $$\triangle U (x)= -\langle C_{1},\lambda_{x}\rangle
=-\dot{X}^{\overline{\beta}}_{x}, x\in
  D ,$$
  Thus we finish the proof of the theorem.  \hfill{$\Box$} \\
\noindent  \\[4mm]
\noindent {\bbb 5\quad The stochastic linear heat   equation  driven by d-parameter fractional L\'{e}vy noise }\\[1mm]
In this section, we consider the  linear stochastic heat   equation
driven by d-parameter fractional L\'{e}vy noise:
$$
\begin{cases}
  \frac{ \partial}{\partial t} U (t, x)=\frac{1}{2} \triangle U (t, x)+\dot{X}^{ \beta_{0},\beta_{1},\ldots,\beta_{d}}_{t,x}, x\in
  D,
\\
U (0, x)=0,t>0,x\in
  D,
\\
 U (t,x)=0, x\in\partial D.
\end{cases}\eqno(5.1)
$$
where $0<\beta_{k}<\frac{1}{2}, k=0,1,\ldots,d$,
$\triangle=\sum_{k=1}^{d}\frac{\partial^{^{2}}}{\partial
x_{k}^{2^{}}}$ is the Laplace operator in $\mathbb{R}^{d}$,
$D\subset \mathbb{R}_{+}^{d}$ is a given domain with regular
boundary and $\dot{X}^{ \beta_{0},\beta_{1},\ldots,\beta_{d}}_{t,x}$
is the  $d+1$-parameter fractional L\'{e}vy noise.
\par
Based on the corresponding solution in the deterministic case, we
guess that
 $$
U(t,x)=\int_{0}^{t}\int_{D}G_{t-s}(x,y)\dot{X}^{
\beta_{0},\beta_{1},\ldots,\beta_{d}}_{s,y}dyds,\eqno(5.2)
$$
where G is the Green function of the heat operator. In fact, we can prove that $U$ is the unique strong solution.\\
\textbf{Theorem 5.1}\ \  The stochastic heat equation (5.1) has a
unique strong  solution in $U:[0,\infty)\times
D\longrightarrow(\mathcal{S})^{-1}$. The solution is
$$
U(t,x)=\int_{0}^{t}\int_{D}G_{t-s}(x,y)\dot{X}^{
\beta_{0},\beta_{1},\ldots,\beta_{d}}_{s,y}dyds,
$$
where G is the Green function of the heat operator $\frac{
\partial}{\partial t}-\frac{1}{2} \triangle $, and (5.2) belongs to
 $C^{1,2}([0,\infty)\times D,(\mathcal{S})^{-1})\cap C ([0,\infty)\times \overline{D},(\mathcal{S})^{-1})$.  \\
\textbf{Proof}:  We first prove that
$U(t,x)\in(\mathcal{S})^{-1}_{-p}$ for all $p>\max\{1,p_{0}\}$.
$$
U(t,x)=\int_{0}^{t}\int_{D}G_{t-s}(x,y)\dot{X}^{
\beta_{0},\beta_{1},\ldots,\beta_{d}}_{s,y}dyds=\int_{0}^{t}\int_{D}G_{t-s}(x,y)\langle
C_{1}, \lambda_{s,y}\rangle dyds
$$
Then by (3.16) ,  for $p>\max\{1,p_{0}\}$,
 $$
\|U(t,x)\|_{-1,-p}
 \leq \int_{0}^{t}\int_{D}\|\langle C_{1},\lambda_{s,y}\rangle\|_{-1,-p} G_{t-s}(x,y) dy <+\infty,
  $$
   that is, $U(t,x)\in(\mathcal{S})^{-1}_{-p}$ for
all $p>\max\{1,p_{0}\}$, for all $t, x$ and
$$
U(t,x)=\langle C_{1},
\int_{0}^{t}\int_{D}G_{t-s}(x,y)\lambda_{s,y}dyds\rangle.\eqno
$$
In fact, the estimate also shows that $U(t,x)$ is uniformly continuous
function from $[0,T]\times \overline{D}$ into $(\mathcal{S})^{-1}$
for any $T<\infty$. Moreover, by the properties of the operator
$G_{t-s}(x,y)$, we get from (5.2) that
$$\aligned
&\frac{\partial}{\partial t}U(t,x)-\frac{1}{2} \triangle U(t,x)\\
&=\dot{X}^{
\beta_{0},\beta_{1},\ldots,\beta_{d}}_{t,x}+\int_{0}^{t}\int_{D}(\frac{
\partial}{\partial t}-\frac{1}{2} \triangle) G_{t-s}(x,y)dyds\\
&=\dot{X}^{ \beta_{0},\beta_{1},\ldots,\beta_{d}}_{t,x}
\endaligned $$
So, $U(t,x)$ satisfies (5.1).      \hfill{$\Box$}
 \par Moreover,we
can prove that under some condition, the solution $U(t,x)$ of
(5.1) is $L^{2}$-integrable.\\
\textbf{Theorem 5.2}\ \ If $
2\beta_{0}+\sum_{i=1}^{d}\beta_{i}+1>\frac{d}{2} $,
then $U(t,x)\in L^{2}(\Omega)$ for all $t\geq 0$, $x\in \overline{D}$.\\
\textbf{Proof}:  From \cite{Hu}, we know that $G$ is smooth in
$(0,\infty)\times D$ and that in $(0,\infty)\times D$,
$$\aligned
& |G_{u}(x,y)|\sim u^{-\frac{d}{2}}\exp(\frac{|x-y|^{2}}{\delta u})\\
& |\frac{\partial G_{u}(x,y)}{\partial y_{i}}|\sim
u^{-\frac{d}{2}-1}|x_{i}-y_{i}|\exp(\frac{|x-y|^{2}}{\delta
 u}),\endaligned
$$
the notion $X\sim Y$ in $(0,\infty)\times D$  means that
$\frac{1}{C}X\leq Y \leq CX$ for some positive constant $C<\infty$
depending only on $D$. By this result, We use the similar proof of Theorem 8.4.1 of \cite{Hu} to verify  the condition
for $U(t,x)\in L^{2}$ for all $t\geq 0$, $x\in \overline{D}$.
$$
\aligned   \mathbb{E}(U(t,x))^{2} & =\mathbb{E}
(\int_{0}^{t}\int_{D}G_{t-s}(x,y)\dot{X}^{
\beta_{0},\beta_{1},\ldots,\beta_{d}}_{s,y}dyds )^{2} \\
&=\mathbb{E} (\int_{0}^{t}\int_{D}(I^{
\beta_{0},\beta_{1},\ldots,\beta_{d}}_{--\ldots-}G_{t-\cdot}(x,\cdot))^{2}dyds )\\
&\sim\int_{0}^{t}\int_{0}^{t}\int_{D}\int_{D}|G_{t-s}(x,y)||G_{t-r}(x,y)||r-s|^{2 \beta_{0} -1}\\
&\quad \quad \prod_{i=1}^{d}|y_{i}-z_{i}|^{2\beta_{i}-1}dy_{1}\ldots
dy_{d}dz_{1}\ldots dz_{d}ds dr \\
&\sim\int_{0}^{t}\int_{0}^{t}\int_{D}\int_{D} (t-r) ^{-\frac{d}{2}}
e^{-\frac{|x-y|^{2}}{\delta (t-r)}} (t-s) ^{-\frac{d}{2}}e^{-\frac{|x-y|^{2}}{\delta(t-s)}}|r-s|^{2 \beta_{0} -1}\\
&\quad \quad \prod_{i=1}^{d}|y_{i}-z_{i}|^{2\beta_{i}-1}dy_{1}\ldots
dy_{d}dz_{1}\ldots dz_{d}ds dr \\
&=\int_{0}^{t}\int_{0}^{t}\int_{D}\int_{D}  r^{-\frac{d}{2}}  s^{-\frac{d}{2}}e^ {-\frac{|x-y|^{2}}{\delta r }}e^ {-\frac{|x-y|^{2}}{\delta s}}|r-s|^{2 \beta_{0} -1}\\
&\quad \quad \prod_{i=1}^{d}|y_{i}-z_{i}|^{2\beta_{i}-1}dy_{1}\ldots
dy_{d}dz_{1}\ldots dz_{d}ds dr. \\
\endaligned \eqno(5.3)
$$
By the inequality (2.1) of \cite{memin} ,
$$
\int_{\mathbb{R}}\int_{\mathbb{R}}|f(x)||g(y)||x-y|^{2 \beta -1}\leq
C
\|f\|_{\frac{1}{\beta+\frac{1}{2}}}\|g\|_{\frac{1}{\beta+\frac{1}{2}}},
0<\beta<\frac{1}{2},\eqno(5.4)
$$
where $C$ is a positive constant, we  have
$$\aligned
&\prod_{i=1}^{d}\int_{-\frac{1}{2}R}^{\frac{1}{2}R}\int_{-\frac{1}{2}R}^{\frac{1}{2}R}e^
{-\frac{|x_{i}-y_{i}|^{2}}{\delta  r }}e^
{-\frac{|x_{i}-y_{i}|^{2}}{\delta s}}(y_{i}-z_{i})^{2\beta_{i}-1}dy_{i}dz_{i} \\
&\leq \prod_{i=1}^{d}[\int_{-\frac{1}{2}R}^{\frac{1}{2}R}e^
{-\frac{|x_{i}-y_{i}|^{2}}{(\beta_{i}+\frac{1}{2})\delta  r
}}dy_{i}]^{\beta_{i}+\frac{1}{2}}[\int_{-\frac{1}{2}R}^{\frac{1}{2}R}e^
{-\frac{|x_{i}-z_{i}|^{2}}{(\beta_{i}+\frac{1}{2})\delta  r }}dz_{i}
]^{\beta_{i}+\frac{1}{2}}\\
&\sim (rs)^ {\frac{1}{2} \sum_{i=1}^{d}(\beta_{i}+\frac{1}{2})},
\endaligned \eqno(5.5 )
$$
where $R$ is some constant such that $D\subset
[-\frac{1}{2}R,\frac{1}{2}R]^{d}$. Substituting (5.5) into (5.3), we
have
$$
\mathbb{E}(U(t,x))^{2}\leq C\int_{0}^{t}\int_{0}^{t}(rs)^
{\frac{1}{2} \sum_{i=1}^{d}(\beta_{i}+\frac{1}{2})}|r-s|^{2
\beta_{0} -1}<\infty
$$
if $ 2\beta_{0}+\sum_{i=1}^{d}\beta_{i}+1>\frac{d}{2} $. Thus we
complete the proof.\hfill $\Box$

\noindent  \\[4mm]
\noindent {\bbb 6\quad The quasi-linear stochastic fractional  heat   equation
driven by {\textcolor{red}{$d$-parameter}} fractional L\'{e}vy noise }\\[1mm]
In this section, we consider the following quasi-linear equation
driven by {\textcolor{red}{$d$-parameter}} fractional L\'{e}vy
noise:
$$
\begin{cases}
  \frac{ \partial}{\partial t} U (t, x)=\frac{1}{2} \triangle U (t, x)+f(U (t, x))+\dot{X}^{ \beta_{0},\beta_{1},\ldots,\beta_{d}}_{t,x}, t>0, x\in
 \mathbb{R}^{d},
\\
U (0, x)=U _{0}(x),x\in
  \mathbb{R}^{d},
\end{cases}\eqno(6.1)
$$
where $U _{0}(x)$ is a given bounded deterministic function on
$\mathbb{R}^{d}$, $f:\mathbb{R}\rightarrow \mathbb{R}$ is a function
satisfying
$$
|f(x)-f(y)|\leq L|x-y|,\forall x, y\in \mathbb{R},\eqno(6.2)
$$
$$
|f(x) |\leq C(1+|x|),\forall x \in\mathbb{R}.\eqno(6.3)
$$
$ U (t, x)$ solves (6.1) if and only if it solves the following
integral equation:
$$\aligned
U (t, x)&=\int_{\mathbb{R}^{d}}U
_{0}(y)G_{t}(x,y)dy+\int_{0}^{t}\int_{\mathbb{R}^{d}}f(U (s,
y))G_{t-s}(x,y)dyds\\
&+\int_{0}^{t}\int_{\mathbb{R}^{d}}G_{t-s}(x,y)dX^{
\beta_{0},\beta_{1},\ldots,\beta_{d}}_{t,x},\endaligned\eqno(6.4)
$$
where
$$
G_{t-s}(x,y)=(t-s)^{-\frac{d}{2}} e^{-\frac{|x-y|^{2}}{2 (t-s)}},
s<t,x\in
  \mathbb{R}^{d},
$$
is the Green function for the heat operator $ \frac{
\partial}{\partial t}-\frac{1}{2} \triangle$.\\
\textbf{Theorem 6.1}\ \ If
$$
\beta_{i}>\frac{1}{2}-\frac{1}{d},i=  1, \ldots,d, \eqno(6.5)
$$
then there exist a unique solution $U(t,x)$ for ( 6.1) such that $U(t,x)\in L^{2}(\Omega)$ for all $t\geq 0$, $x\in \overline{D}$.\\
\textbf{Proof}: Define
$$
V(t,x)=\int_{0}^{t}\int_{\mathbb{R}^{d}}G_{t-s}(x,y)dX^{
\beta_{0},\beta_{1},\ldots,\beta_{d}}_{t,x}.
$$
Since (6.5) holds, by the similar arguments of Theorem 5.2, we can
prove that $ V(t,x) \in L^{2}(\Omega)$ for all $t\geq 0$, $x\in
\overline{D}$, so $ V(t,x)$ exists as an ordinary random field. The
existence of the solution now follows Picard iteration. Define
 $$
 U_{0}(t,x)=U _{0}(x)
 $$
and iteratively
$$\aligned
U_{j+1} (t, x)&=\int_{\mathbb{R}^{d}}U
_{0}(y)G_{t}(x,y)dy+\int_{0}^{t}\int_{\mathbb{R}^{d}}f(U_{j} (s,
y))G_{t-s}(x,y)dyds\\
&+\int_{0}^{t}\int_{\mathbb{R}^{d}}G_{t-s}(x,y)dX^{
\beta_{0},\beta_{1},\ldots,\beta_{d}}_{t,x},j=0, 1,
2,\ldots\endaligned\eqno(6.6)
$$
Then by (6.3), $U_{j} (t, x)\in L^{2}(\mathbb{P})$ for all $j$. And
by (6.2),
$$\aligned
&\mathbb{E}|U_{j+1} (t, x)-U_{j} (t,x)|^{2}\\
&=\mathbb{E}|\int_{0}^{t}\int_{\mathbb{R}^{d}}(f(U_{j} (s,y))-f(U_{j-1} (s, y)))G_{t-s}(x,y)dyds |^{2}\\
& \leq L\mathbb{E}[\int_{0}^{t}\int_{\mathbb{R}^{d}}|U_{j} (s,y)-U_{j-1} (s, y)|G_{t-s}(x,y)dyds] ^{2}\\
& \leq L\int_{0}^{t}\int_{\mathbb{R}^{d}}G_{t-s}(x,y)dyds
\int_{0}^{t}\int_{\mathbb{R}^{d}}\mathbb{E}|U_{j} (s,y)-U_{j-1} (s,
y)|^{2}G_{t-s}(x,y)dyds\\
& \leq C_{T} \int_{0}^{t} \sup_{y}\mathbb{E}|U_{j} (s,y)-U_{j-1} (s,
y)|^{2} ds\\
& \leq C_{T}^{j} \int_{0}^{t} \int_{0}^{s_{1}}\ldots
\int_{0}^{s_{j-1}} \sup_{y}\mathbb{E}|U_{1} (s,y)-U_{0} (s,
y)|^{2}ds_{j-1}\ldots ds_{1}ds\\
& \leq A_{T}C_{T}^{j}\frac{T^{j}}{j!}
\endaligned
$$
for some constants $A_{T},C_{T}$. It follows that the sequence
$\{U_{j} (t, x) \}_{j=0}^{\infty}$ of random fields converges in
$L^{2}(\mathbb{P})$ to a random field $U(t,x)$. Letting
$j\rightarrow\infty$ in (6.6), we see that $U(t,x)$ is a solution of
(6.1). The uniqueness follows by the Gronwall$^{,}$s
inequality.\hfill $\Box$\\
 \\[4mm]

 \noindent{\bbb{References}}

\end{document}